\documentclass{article}
\usepackage{spconf, amsmath, graphicx}

\usepackage[utf8]{inputenc} 
\usepackage[T1]{fontenc}    
\usepackage{hyperref}       
\usepackage{url}            
\usepackage{booktabs}       
\usepackage{amsfonts}       
\usepackage{nicefrac}       
\usepackage{microtype}      
\usepackage{multirow}
\usepackage{cite}

\usepackage{ amssymb }
\usepackage{ amsmath }
\usepackage{ amsfonts}
\usepackage{ dsfont }
\usepackage{ amsthm}
\usepackage{algpseudocode}
\usepackage[ruled, linesnumbered]{algorithm2e}
\usepackage{graphicx}
\usepackage{bbm}
\usepackage{ enumerate}
\usepackage{bm}
\usepackage{tikz}
\usepackage{stmaryrd}

\newcommand{\C}{\mathbb{C}}

\newcommand{\bra}{\langle}
\newcommand{\ket}{\rangle}

\newcommand{\bPhi}{\boldsymbol \Phi}
\newcommand{\bvarphi}{\boldsymbol \varphi}

\newcommand{\bgamma}{\boldsymbol \gamma}

\newcommand{\beps}{{\boldsymbol \epsilon}}

\newcommand{\eqdef}{\overset{\text{\tiny{def}}}{=}}

\newcommand{\act}{\mathcal{A}} 
\newcommand{\supp}{\mathcal{S}} 

\newcommand{\bi}{\mathcal{H}} 

\newtheorem{thm}{Theorem}
\newtheorem{defi}[thm]{Definition}

\newtheorem{cor}[thm]{Corollary}

\numberwithin{thm}{section}
\numberwithin{equation}{section}
\numberwithin{figure}{section}

\DeclareMathOperator\range{col}

\title{On the Optimality of Backward Regression:\\
 Sparse Recovery and Subset Selection
}
\name{Sebastian Ament \qquad and \qquad Carla Gomes 
\thanks{This research was supported by NSF awards CCF-1522054 (Expeditions in computing) and  AFOSR Multidisciplinary University Research Initiatives (MURI) Program FA9550-18-1-0136, ARO award W911NF-17-1-0187 for our compute cluster, US DOE Award No. DE-SC0020383, and an award from the Toyota Research Institute.}}
\address{Cornell University \\
Department of Computer Science \\ 
Ithaca, NY 14850}

\begin{document}

\maketitle
\begin{abstract}
Sparse recovery
and subset selection
are fundamental problems in varied communities,
including signal processing, statistics and machine learning.
Herein, we focus on an important greedy algorithm for these problems:
Backward Stepwise Regression.
We present novel guarantees for the algorithm,
propose an efficient, numerically stable implementation,
and 
put forth Stepwise Regression with Replacement (SRR),
a new family of two-stage algorithms that employs both forward and backward steps
for compressed sensing problems.
Prior work on the backward algorithm
has proven its optimality for the subset selection problem,
provided the residual associated with the optimal solution is small enough.
However, the existing bounds on the residual magnitude are NP-hard to compute.
In contrast, our main theoretical result includes a bound that can be computed in polynomial time, depends chiefly on the smallest singular value of the matrix,
and also extends to the method of magnitude pruning.
In addition, we report numerical experiments 
highlighting crucial differences between 
forward and backward greedy algorithms 
and compare SRR against popular two-stage algorithms 
for compressed sensing. 
Remarkably, SRR algorithms generally maintain good sparse recovery performance 
on coherent dictionaries.
Further, a particular SRR algorithm has an edge over Subspace Pursuit.
\end{abstract}


\section{Introduction}
\label{sec:intro}

Sparse signal recovery and subset selection
are problems with important applications in
medicine \cite{rao2012alzheimer, rao2013fetal}, 
materials science \cite{ghiringhelli2017learning},
and engineering \cite{chan2008single}.
Given a set of {\it atoms} $\mathcal{D} \eqdef \{\bvarphi_i\}$,
also referred to as a dictionary,
both problems are about identifying a sparse subset of $\mathcal{D}$ 
that most accurately model an observation $\bold y$, but the assumptions vary.
Formally, the central problem is the solution of the following expression with a desired sparsity:
\begin{equation}
\label{eq:problem}
\min_{|\act| = k} \min_{\bold x} \|\bold y - \sum_{i \in \act} x_i \bvarphi_i \|_2
\end{equation}
If $\bold y = \bPhi \bold x + \beps$ 
where $\bPhi \eqdef [ \bvarphi_1, ..., \bvarphi_m]$,
 $\bold x$ has at most $k$ non-zero entries,
 and $\beps$ is a small perturbation,
the problem is referred to as sparse recovery.
On the other hand, if $\bold y$ is an arbitrary vector,
not necessarily associated with a sparse $\bold x$,
the problem is called sparse approximation
or subset selection.
Notably, the recent work of \cite{rao2020ssgd}
generalizes the framework of sparse recovery
for the sparsification of neural networks
with promising results.

Because of the importance and ubiquity of the problems, 
many approaches have been proposed in the literature 
\cite{candes2006robust, greedisgood, ssp, cosamp, wipf2008, gradientpursuit2008, blumensaththres2009, he2017bayesian, tipping2001sparse}.
Two popular greedy algorithms are Forward and Backward Regression,
which have been proposed repeatedly in different fields,
 leading to a perplexing number of names 
that all refer to the same two algorithms.  

In particular, the forward algorithm is also known as 
Forward Regression, Forward Selection \cite{miller2002subset}, 
Optimized Orthogonal Matching Pursuit (OOMP) \cite{rebollo2002oomp},
Order-recursive Matching Pursuit (ORMP) \cite{cotter1999ormp},
and Orthogonal Least-Squares (OLS) \cite{blumensath2007difference, chen1989ols, chen1991ols}.
The backward algorithm is also known as 
Backward Regression, Backward Elimination \cite{miller2002subset},
Backward Optimized Orthogonal Matching Pursuit (BOOMP) \cite{andrle2004boomp},
and the backward greedy algorithm \cite{backwardoptimality}.
The two algorithms add atoms to, respectively delete atoms from, an {\it active set} $\act$, based on the greedy heuristics
\begin{equation}
\label{eq:stepwise}
\arg \min_{i \not \in \act} \| \bold r_{\act \cup i} \|_2, 
\qquad 
\text{and}
\qquad
\arg \min_{i \in \act} \| \bold r_{\act \backslash i} \|_2,
\end{equation}
where 
$\bold r_\act = (\bold I - \bold \bPhi_\act \bPhi_\act^+) \bold y$
is the least-squares residual associated with
the subset of atoms indexed by $\act$.
Among all additions and deletions of a single atom,
these heuristics leads to the minimum residual in the following iteration.
By comparison, the well-known Orthogonal Matching Pursuit (OMP) algorithm
adds atoms based on the rule
\begin{equation}
\label{eq:omp}
\arg \max_{i \not \in \act} |\bra \bvarphi_i, \bold r_\act \ket|,
\end{equation} 
which is equivalent to the largest component in the gradient 
of the least-squares objective at the current residual,
 $\arg \max_{i \not \in \act} \left| \partial_{x_i} \| \bold r_\act \| \right|$.
In this sense, the stepwise algorithms "look ahead" further than OMP.

Herein, our primary focus is the backward algorithm,
for which we propose an efficient and numerically stable implementation and
derive novel theoretical insights which also extend to magnitude pruning.
An important limitation of the algorithm is that it only works 
with matrices $\bPhi$ that have full column rank.
Perhaps because of this limitation and 
the high dimensionality of modern datasets,
 recent focus has been primarily on forward algorithms.
To overcome this limitation and 
take advantage of the strong theoretical guarantees of the backward algorithm,
we further propose Stepwise Regression with Replacement,
a new family of two-stage algorithms that uses
both forward and backward steps as building blocks
to solve challenging compressed sensing problems
efficiently.

\section{PRIOR WORK}

We begin by reviewing the most relevant results on the optimality of the backward algorithm
and make connections to other work whenever it becomes relevant.
A central concept to the existing theory
is the bisector of two linear spaces.

\begin{defi}[Bisector]
Let $A$ and $B$ be two subspaces of a normed linear space $L$ 
with metric $d$.
Their bisector is
\[
\bi(A, B) \eqdef \{ \bold x \in L | d(\bold x, A) = d(\bold x, B) \},
\]
where 
 $d(\bold x, S) \eqdef \min_{\bold s \in S} d(\bold x, \bold s)$ for any $S \subset L$.
\end{defi}
The proof of the following result, due to Covreur and Bresler \cite{backwardoptimality},
views certain bisectors 
as decision boundaries of the backward algorithm,
and shows that, as long as the boundaries cannot be crossed due to the perturbation,
the algorithm succeeds in recovering a sparse signal.

\begin{thm}[Covreur and Bresler \cite{backwardoptimality}]
\label{thm:backwardoptimality}
Backward Regression with input $\bPhi \bold x + \beps$ 
recovers the support $\supp$ of 
a $k$-sparse $m$-dimensional vector $\bold x$ 
in $m-k$ iterations provided
\[
\min_{k < r \leq n} \min_{|\act| = r} \min_{i \in \act, j \not \in \act} 
d(\bold y,  \bi(\range(\bPhi_{\act \backslash i}), \range(\bPhi_{\act \backslash j})))
> \| \beps \|_2,
\]
where $\range(\bPhi)$ denotes the column space of $\bPhi$.
\end{thm}
The authors postulate that the bound in Theorem \ref{thm:backwardoptimality} is NP-hard to compute and therefore of limited practical use.
Further, Covreur and Bresler proved that
the backward algorithm is not only capable of recovering the support of an {\it exactly} sparse vector, but that it solves the subset selection problem optimally, 
provided the residual of the optimal solution is small enough.
In particular, for an arbitrary $\bold y$, not generally associated with a sparse $\bold x$, the following result holds.
\begin{cor}[Covreur and Bresler \cite{backwardoptimality}]
\label{cor:backwardoptimality}
Let $\bold x_k$ be the solution to the subset selection problem \eqref{eq:problem}
with sparsity $k$.
If the residual $\bold r_k \eqdef \bold y - \bPhi \bold x_k$ satisfies the bound in Theorem~\ref{thm:backwardoptimality} in place of $\beps$,
the backward algorithm recovers $\bold x_k$.
\end{cor}

\section{Main Results}

We first propose an efficient and numerically-stable procedure for the evaluation of the column deletion criterion. Subsequently, we present our main theoretical results for the backward algorithm.
In the following, $\bPhi$ is a $n \times m$ matrix.

\subsection{An Efficient and Stable Implementation}
\begin{algorithm}[t]
 \KwData{Matrix $\bPhi \in \C^{n\times m}$, target $\bold y$, desired sparsity $k$}
 \KwResult{Support set $\act$}
 $\act \leftarrow \llbracket 1, m \rrbracket$\\
 $ \bold Q_\act, \bold R_\act \leftarrow \text{qr}(\bPhi_\act)$ \\
 \While{$|\act| > k$}{
 $\bold x_\act \leftarrow \bold R_\act^{-1} \bold Q_\act^* \bold y$ \\
 $\bgamma_\act \leftarrow \text{diag}[(\bold R_\act^* \bold R_\act)^{-1}]$ \\
 $i^* \leftarrow \arg \min_{i \in \act} |x_i|^2 / \gamma_i$ \\
 $\bold Q_{\act\backslash i^*}, \bold R_{\act\backslash i^*} \leftarrow$ remove\_column$(\bold Q_\act, \bold R_\act, i^*)$ \\
   $ \act \leftarrow \act \backslash i^*$ \\
 }
 \caption{Efficient Backward Regression}
 \label{alg:br}
\end{algorithm}

\label{sec:implementation}
Previously, \cite{reeves1999efficient} proposed 
an efficient implementation for the backward algorithm 
based on rank-one updates to the normal equations,
which are notorious for exacerbating ill-conditioning.
Since the algorithm uses recursive updates,
small errors are prone to blow up catastrophically.
To make the backward algorithm usable for larger, potentially ill-conditioned systems,
we require an implementation that is both efficient and numerically stable.
To this end, let
\[
\bgamma_\act \eqdef \text{diag}\left( [\bPhi_\act^* \bPhi_\act ]^{-1} \right).
\]
Using the analysis of \cite{reeves1999efficient}, it can be proven that for $i \in \act$,
\begin{equation}
\label{eq:deltaviagamma}
 \| \bold r_{\act\backslash i} \|_2^2 -  \| \bold r_{\act} \|_2^2 =  | x_i |^2 / \gamma_i.
\end{equation}
 Therefore, we can efficiently evaluate the deletion criterion 
 in \eqref{eq:stepwise}
 by computing $\bgamma$ and the coefficients $\bold x$ 
 corresponding to the least-squares solution given $\bPhi_\act$.
But rather than computing $\bgamma$ and $\bold x$ via the normal equations,
we keep and update a QR factorization of $\bPhi_\act$,
and set $\bgamma \leftarrow \text{diag}[(\bold R^* \bold R)^{-1}]$.
While the computation of $\bgamma$ is potentially still unstable, 
the instability cannot propagate through subsequent iterations.
See Algorithm \ref{alg:br}.

Further, equation \eqref{eq:deltaviagamma} exposes a similarity 
of the backward criterion to {\it magnitude pruning},
a principle used by
Subspace Pursuit \cite{ssp} for its deletion stage,
and popularly applied for the sparsification of deep neural networks
\cite{han2015learning, carreira2018prune, dai2019nest}.
Inspired by this observation,
we also study a modification of Algorithm \ref{alg:br}
where line 6 is replaced with 
\[
\arg \min_{i \in \act} | x_i |.
\]
We refer to the resulting algorithm as Least-Absolute Coefficient Elimination (LACE).
We now highlight the advantages of the proposed implementation.


\begin{figure}[t!]
\includegraphics[width=.238\textwidth]{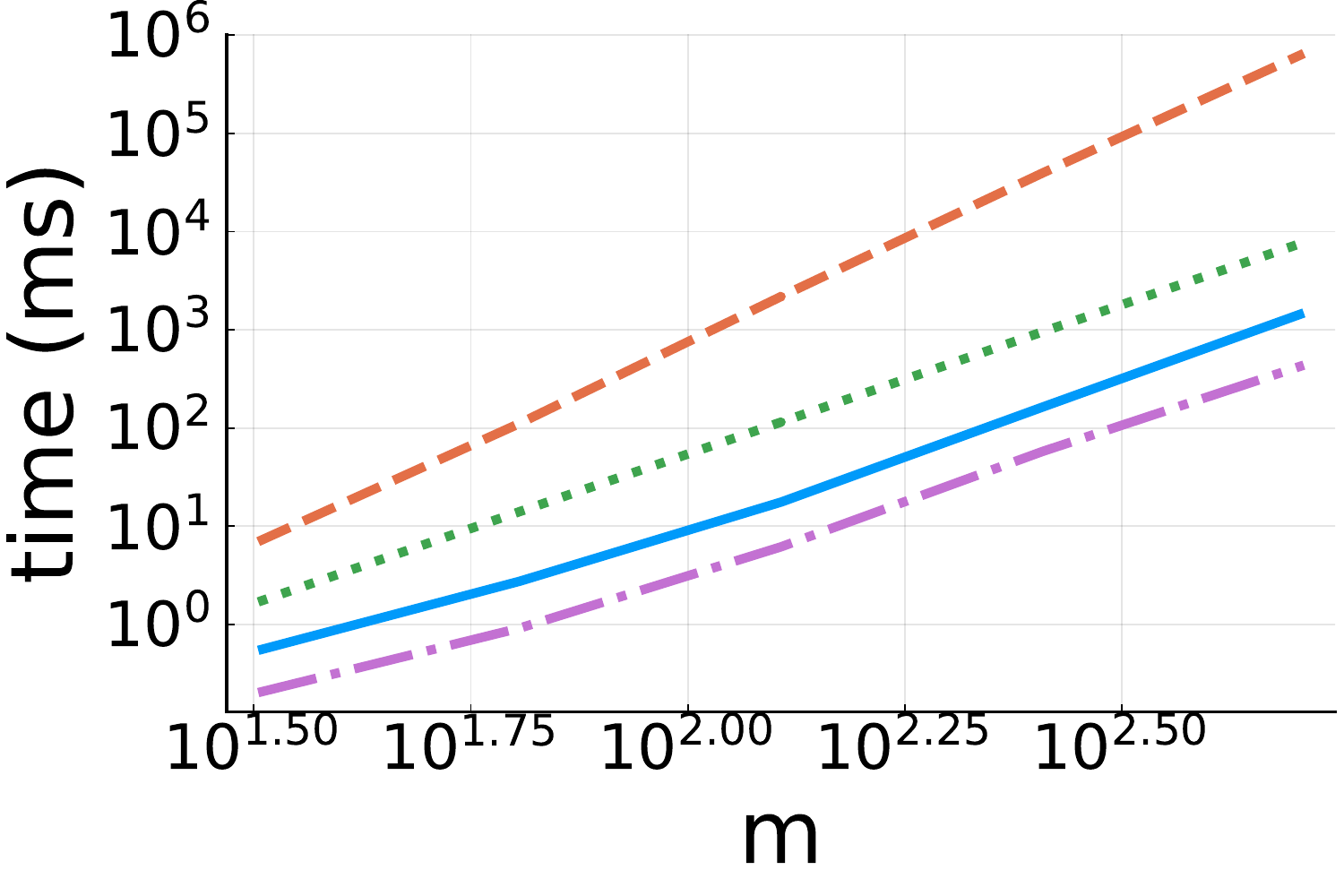}
\includegraphics[width=.238\textwidth]{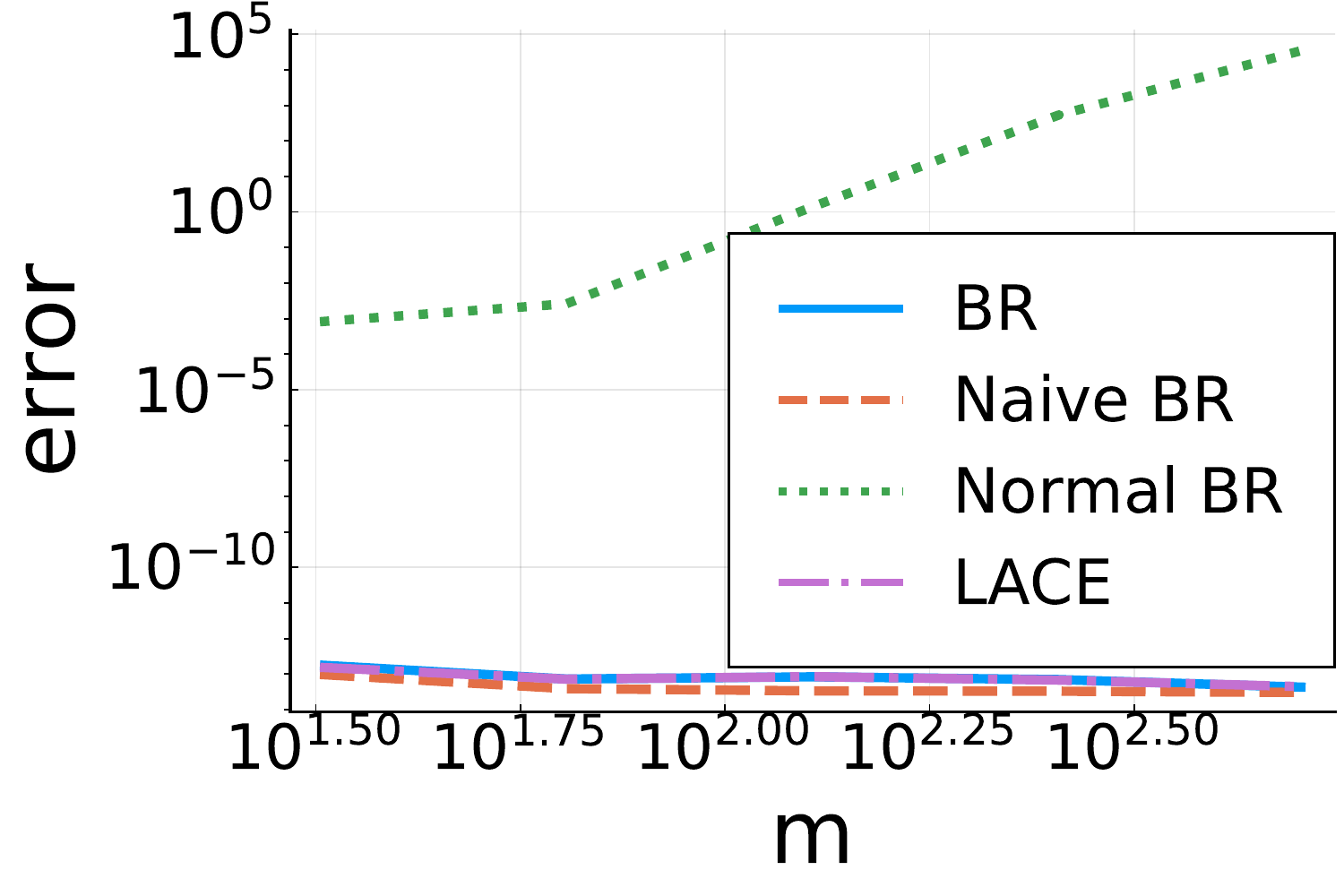}
\caption{Runtime (left) and error (right) as a function of $m$
for ill-conditioned matrices with $k = 16$, no noise, and $n =m$.
}
\label{fig:stability}
\end{figure}

The left of Figure \ref{fig:stability} 
shows representative runtimes of four implemenations: 
Algorithm \ref{alg:br} (BR),
a naive implementation which calculates $\bold r_{\act \backslash i}$
for all $i \in \act$ by 
removing a column of the QR factorization and solving the resulting system (Naive BR),
the efficient implementation of Reeves \cite{reeves1999efficient} (Normal BR),
and LACE.
We fixed the sparsity at $k = 16$ and increased $m$ while keeping the matrices square.
Notably, BR achieves the same scaling as LACE and is only a small constant factor slower,
while the naive implementation is orders of magnitude slower than the other ones.
The right of Figure \ref{fig:stability} 
shows the error of the solution for the same algorithms on randomly generated matrices with a high condition number.
As we did not add any noise, the error should mathematically be zero.
But as predicted, the result of Normal BR  
blows up as the number of iterations of the algorithm grows with $m$.
Consequently, Algorithm \ref{alg:br} is preferable 
for its efficiency and stability.

\subsection{Theoretical Guarantees}
The following
is our main theoretical result, which provides a sparse recovery guarantee for both backward algorithms.

\begin{thm}
\label{thm:berc}
Given $\bPhi$ with full column rank,
Backward Regression and LACE recover the support $\supp$ of a $k$-sparse $m$-dimensional vector $\bold x$ in $m-k$ steps given $\beps$ satisfies
\[
\frac{\sigma_{\min}(\bPhi)}{2} \min_{i \in \supp} |x_i| > \| \beps \|_2,
\]
where $\sigma_{\min}(\bPhi)$ is the smallest singular value of $\bPhi$.
\end{thm}

Theorem \ref{thm:berc} reveals 
the direct dependence of the tolerable noise
on the smallest singular value of the matrix,
something which Theorem \ref{thm:backwardoptimality} obscures
and was only stated heuristically in \cite{backwardoptimality}.
The following corollary is analogous to Corollary \ref{cor:backwardoptimality}
but is based on Theorem \ref{thm:berc} and extends to LACE.
\begin{cor}
\label{cor:berc}
Let $\bold x_k$ be the solution to the subset selection problem \eqref{eq:problem}
with sparsity $k$.
If the associated residual $\bold r_k \eqdef \bold y - \bPhi \bold x_k$ satisfies the bound in Theorem~\ref{thm:berc} in place of $\beps$,
Backward Regression and LACE recover $\bold x_k$.
\end{cor}

While Corollary \ref{cor:berc} cannot guarantee the success of the algorithm a-priori,
it can be used as an efficient post-hoc check of the algorithms' return values.
Alternatively, the algorithms could be stopped adaptively, once 
the residual surpasses the bound to guarantee that the returned value
is optimal among all signals with the same number of non-zero elements.

\begin{figure}[t]
\begin{center}
\includegraphics[width=.215\textwidth]{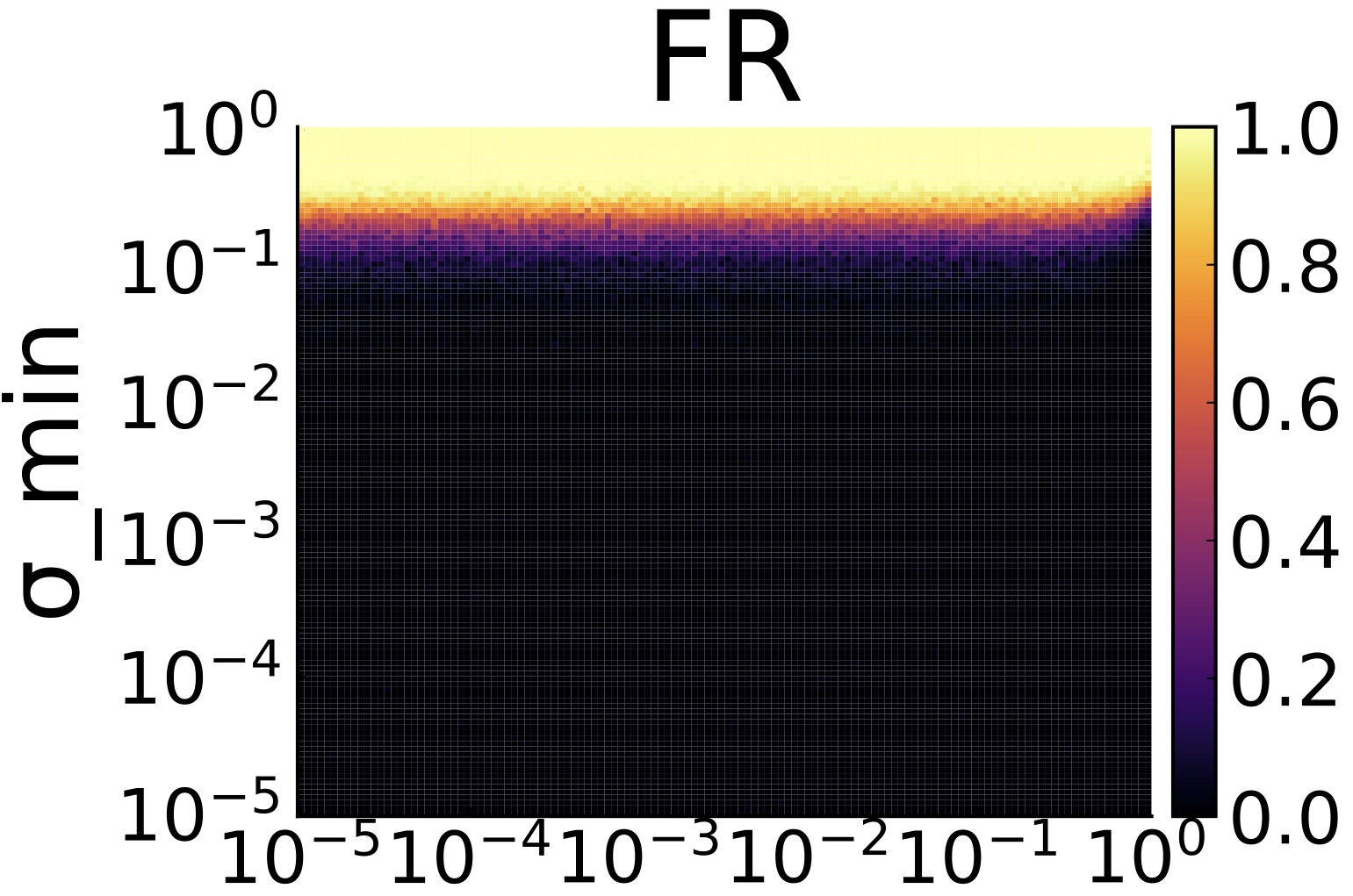}
\includegraphics[width=.205\textwidth]{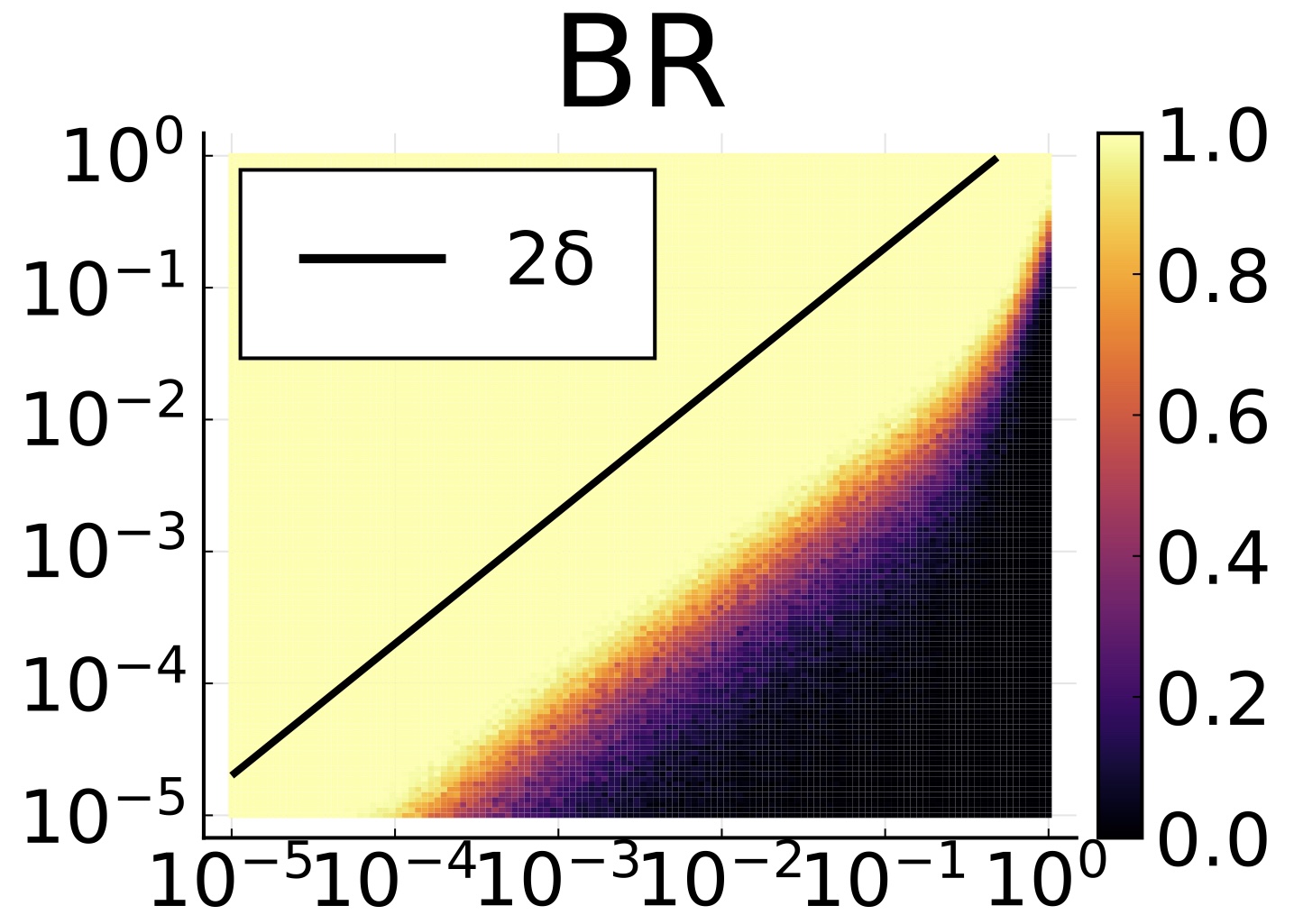}
\includegraphics[width=.215\textwidth]{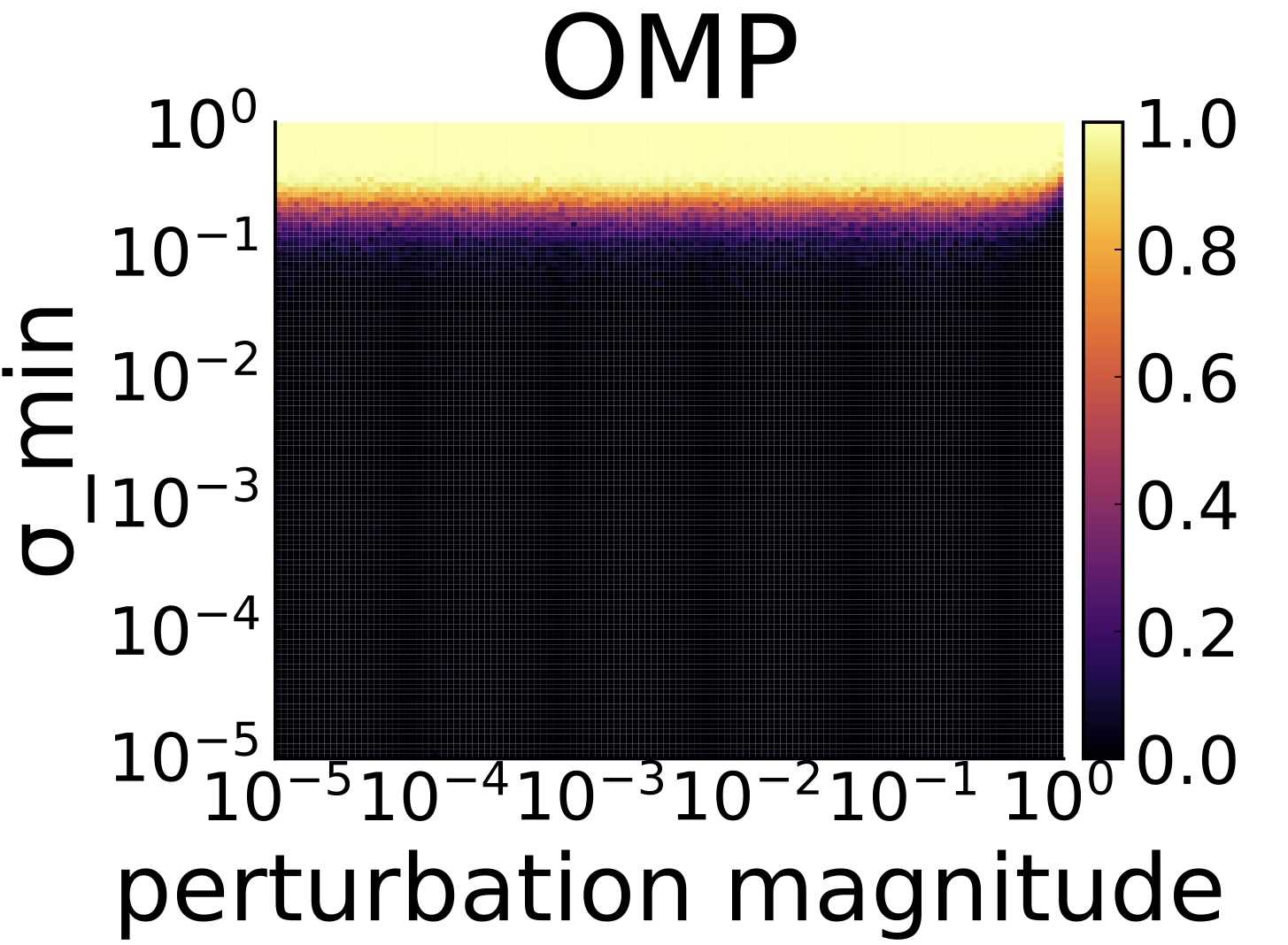}
\includegraphics[width=.205\textwidth]{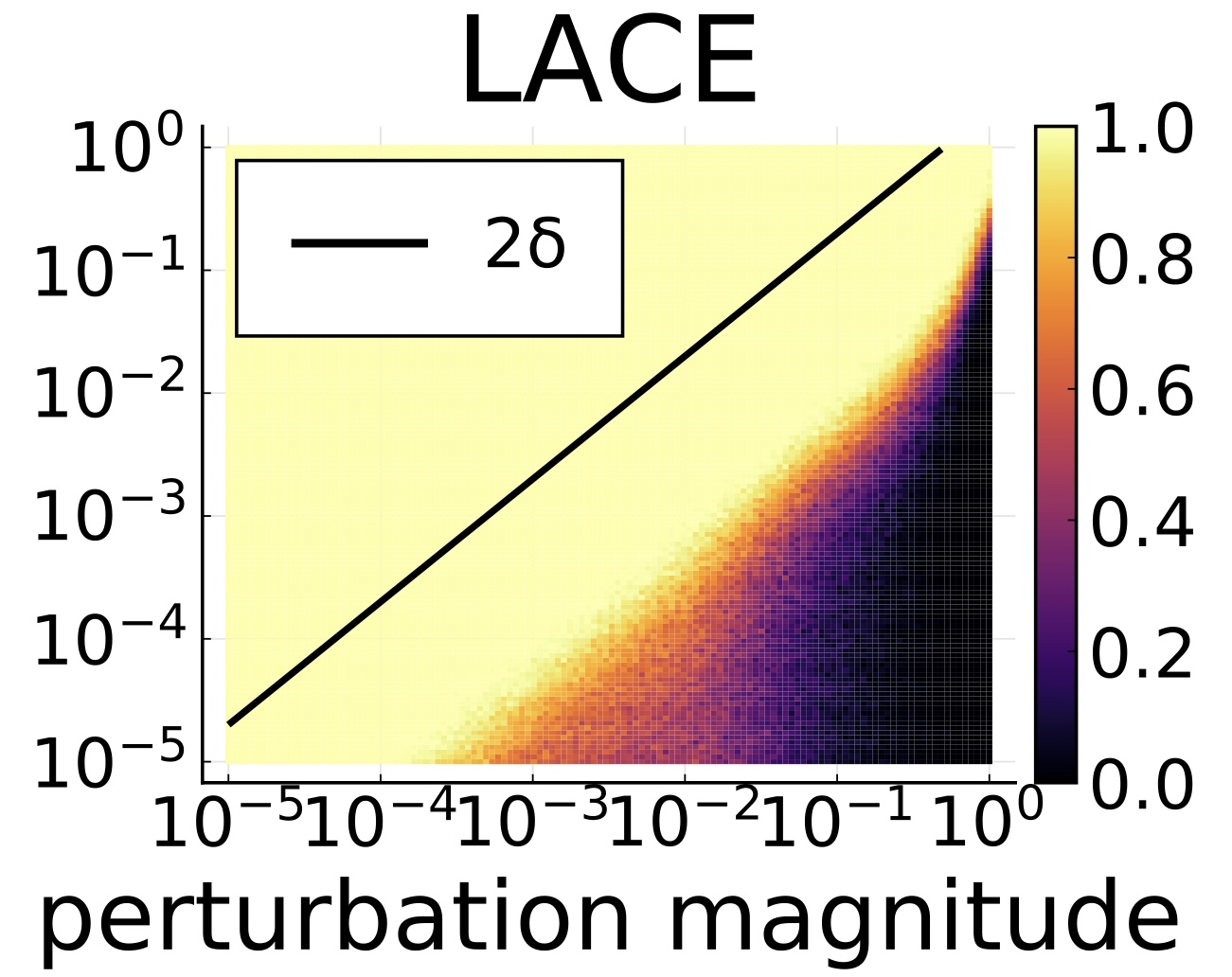}
\caption{Empirical frequency of support recovery as a function of the smallest singular value $\sigma_{\min}$ and perturbation magnitude for 64 by 64 matrices with signal sparsity 32.
}
\end{center}
\label{fig:backward}
\end{figure}

Notably, a necessary condition for OMP and Forward Regression (FR) to recover the support of a sparse signal is $\mu_1(k) < 1/2$, where $\mu_1$ is the Babel function of the matrix $\bPhi$,
whose columns are assumed to be $l_2$-normalized \cite{greedisgood, soussen2013joint}.
This implies 
$\sigma_{\min}(\bPhi_\act) > 1/\sqrt{2} \approx 0.71$
since
for any singular value $\sigma$ of $\bPhi_\act$,
 $|1- \sigma^2 | \leq \mu_1(k)$
for any set $\act$ with cardinality $k$ \cite{greedisgood}.
Therefore, the forward algorithms are not guaranteed to recover a sparse support
for even moderately ill-conditioned systems, regardless of how small the noise is.
This is in stark contrast to Theorem \ref{thm:berc} for the backward algorithms,
which are always guaranteed to return the correct support of a sparse vector
if the noise is small enough.

\subsection{Empirical Evaluation}

Figure \ref{fig:backward} shows the empirical frequency of support recovery 
for Forward Regression (FR), BR, OMP, and LACE
for 64 by 64 matrices
as a function of perturbation magnitude $\delta = \|\beps\|$ and the smallest singular value $\sigma_{\min}$ of $\bPhi$.
For a given $\sigma_{\min}$, 
we define $\bold S_\sigma$ as the matrix with uniformly spaced values between $\sigma_{\min}$ and $1$ on the diagonal, 
and let $\bPhi_\sigma = \bold U \bold S_\sigma \bold V^*$,
where $\bold U, \bold V$ are two orthonormal matrices.
The perturbation vectors are uniformly distributed on the $\delta$-hypersphere,
and $x_i$ for $i \in \supp$ has the Rademacher distribution.
Each cell of the heat-maps is an average of 128 independent experiments.

Remarkably, the recovery performance of both BR and LACE exhibits a scaling that is approximately linear in $\delta$, as
predicted by Theorem \ref{thm:berc},
and surpasses the forward algorithms on 
matrices with small singular values.
While the forward algorithms only succeed on well-conditioned systems,
as suggested by the theory \cite{greedisgood, soussen2013joint},
notably, \cite{das2018approximate}
proved that both FR and OMP are approximation algorithms to the subset selection problem,
using a notion of approximate submodularity 
of the coefficient of determination, $R^2$.
Indeed, additional experiments not reported herein show that the forward algorithms' performance in terms of $R^2$ degrades more gracefully with the smallest singular value, even as the support recovery performance deteriorates sharply.

\section{Stepwise Regression with Replacement}
\label{sec:srr}

\begin{algorithm}[t]
 \KwData{dictionary $\bPhi$, target $\bold y$, sparsity $k$, steps $s$}
 \KwResult{Support set $\act$}
 $\act \leftarrow$ $k$ column indices with largest correlation with $\bold y$ \\
 \While{not converged}{
 \For{$s$ iterations}{
$\act \leftarrow \act \cup \arg \min_{i \not \in \act} \| \bold r_{\act \cup i} \|_2$ 
\Comment{acquisition}
}
 \For{$s$ iterations}{
  $\act \leftarrow \act \backslash \arg \min_{i \not \in \act} \| \bold r_{\act\backslash i} \|_2$ \Comment{deletion}
 }
 }
 \caption{{\small Stepwise Regression with Replacement}}
 \label{alg:srr}
\end{algorithm}

\subsection{Motivation and Algorithm Design}

While the results presented above are powerful
and do not require the same
restrictive assumptions as existing forward and two-stage algorithms  \cite{greedisgood, ssp, cosamp},
the backward algorithms require the matrices to have full column rank.
Thus, they are not immediately applicable for compressed sensing.
To remedy this limitation, we propose Stepwise Regression with Replacement (SRR),
a new family of two-stage algorithms which combines forward and backward steps to replace atoms in an active set. 
See Algorithm \ref{alg:srr} for schematic code of SRR.

Two-stage algorithms generally keep track of an active set with a desired sparsity $k$,
and improve the reconstruction of the observed signal 
by iteratively adding and deleting atoms from the active set.
For example,
both Subspace Pursuit (SP) \cite{ssp}
and Orthogonal Matching Pursuit with Replacement (OMPR) \cite{jain2011orthogonal} 
initialize an active set of size $k$
and subsequently
alternate
adding atoms to the active set based on the largest correlations
with the residual,
and deleting atoms with magnitude pruning. 
SP solves a least-squares problem for the expanded active set before deleting atoms,
while OMPR only takes one gradient step.
Further, SP adds and deletes $k$ atoms,
while OMPR adds and deletes $s$ atoms, 
where $s$ is passed as a parameter.

Stepwise Regression with Replacement 
has a similar basic structure, and also takes an additional "step" parameter $s$.
However, it uses the heuristics in equation \eqref{eq:stepwise} to update the support,
see Algorithm \ref{alg:srr}.
SRR can be implemented efficiently with updates to a QR-factorization,
similar to Algorithm \ref{alg:br}.

\subsection{Empirical Evaluation}
\begin{figure}[t!]
\begin{center}
\includegraphics[width=.22\textwidth]{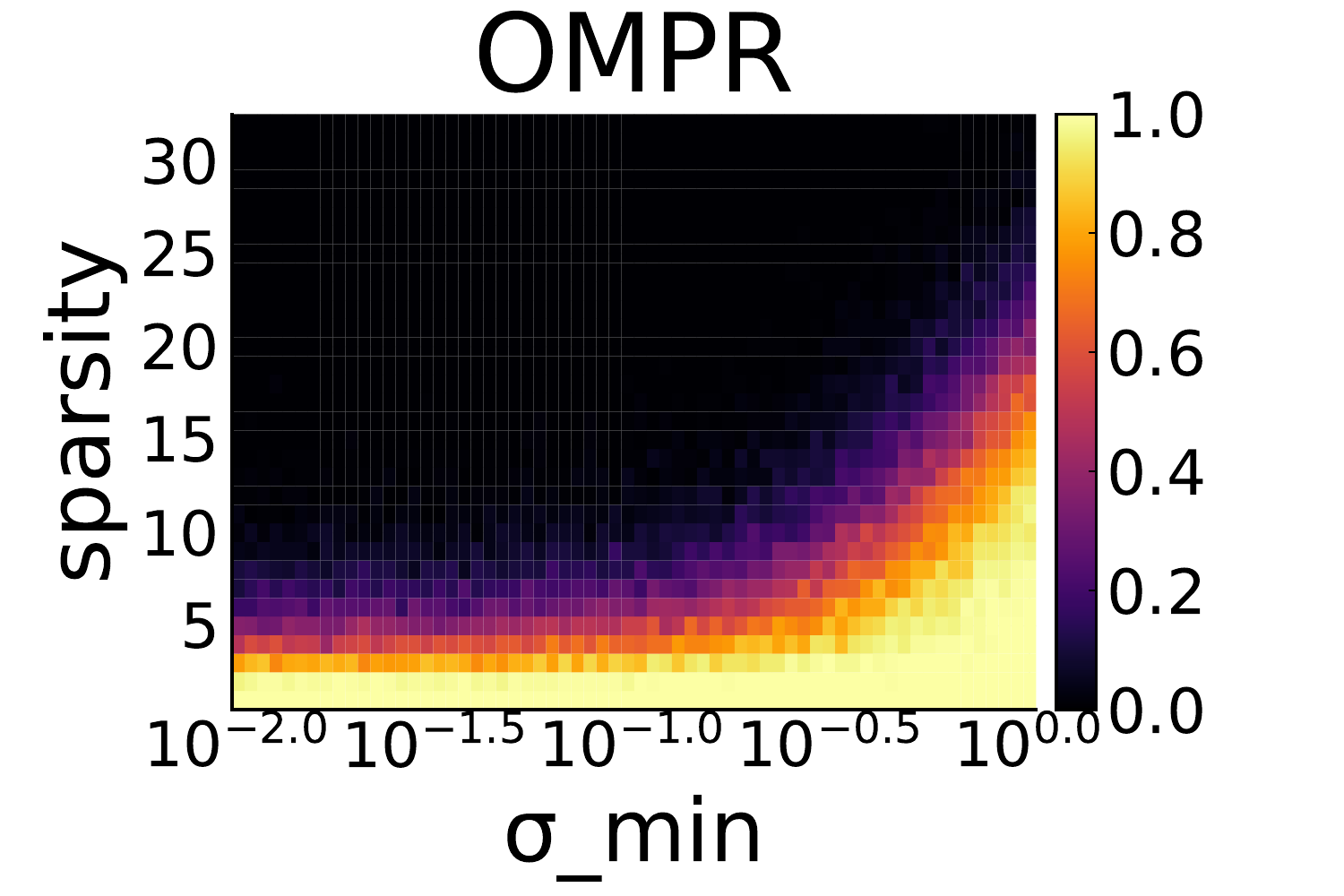}
\includegraphics[width=.205\textwidth]{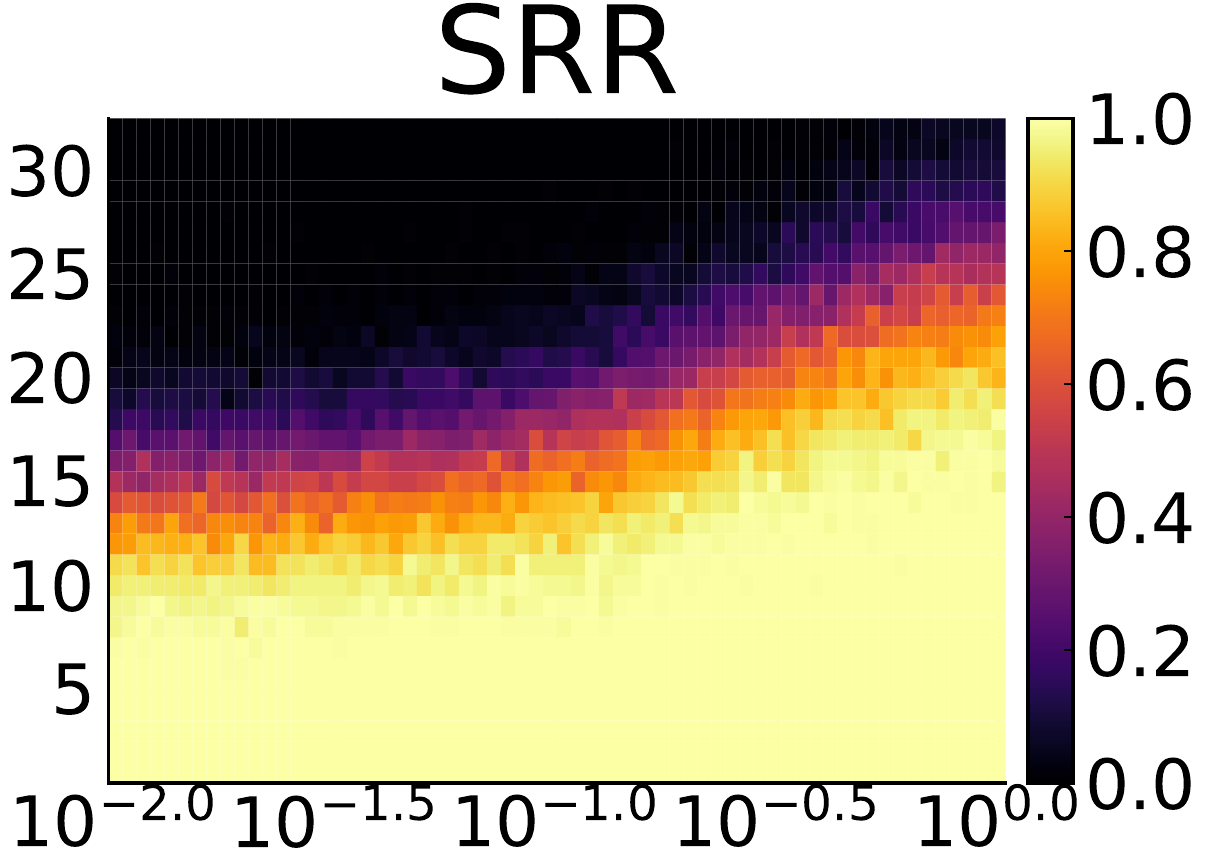}
\includegraphics[width=.22\textwidth]{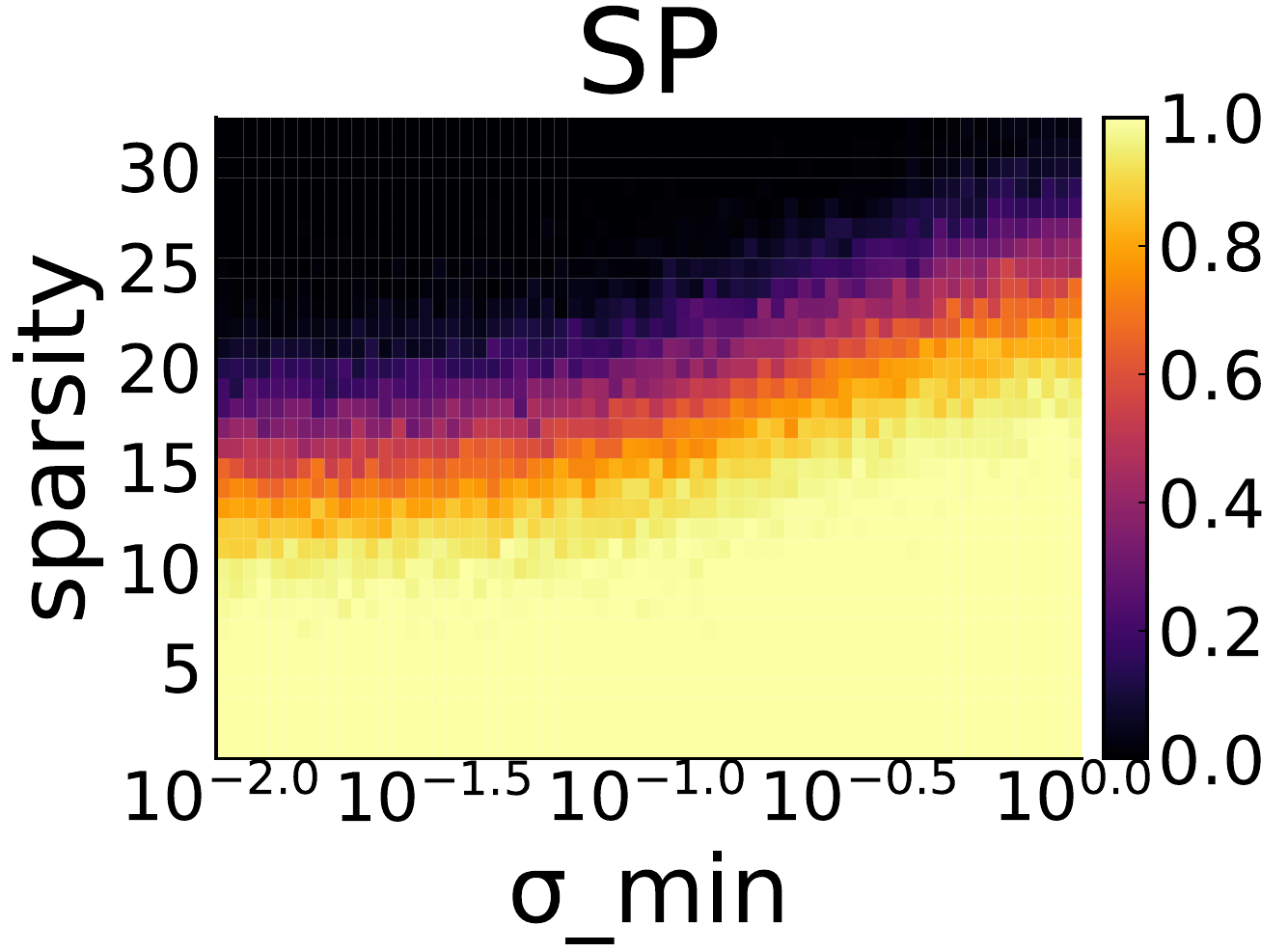}
\includegraphics[width=.205\textwidth]{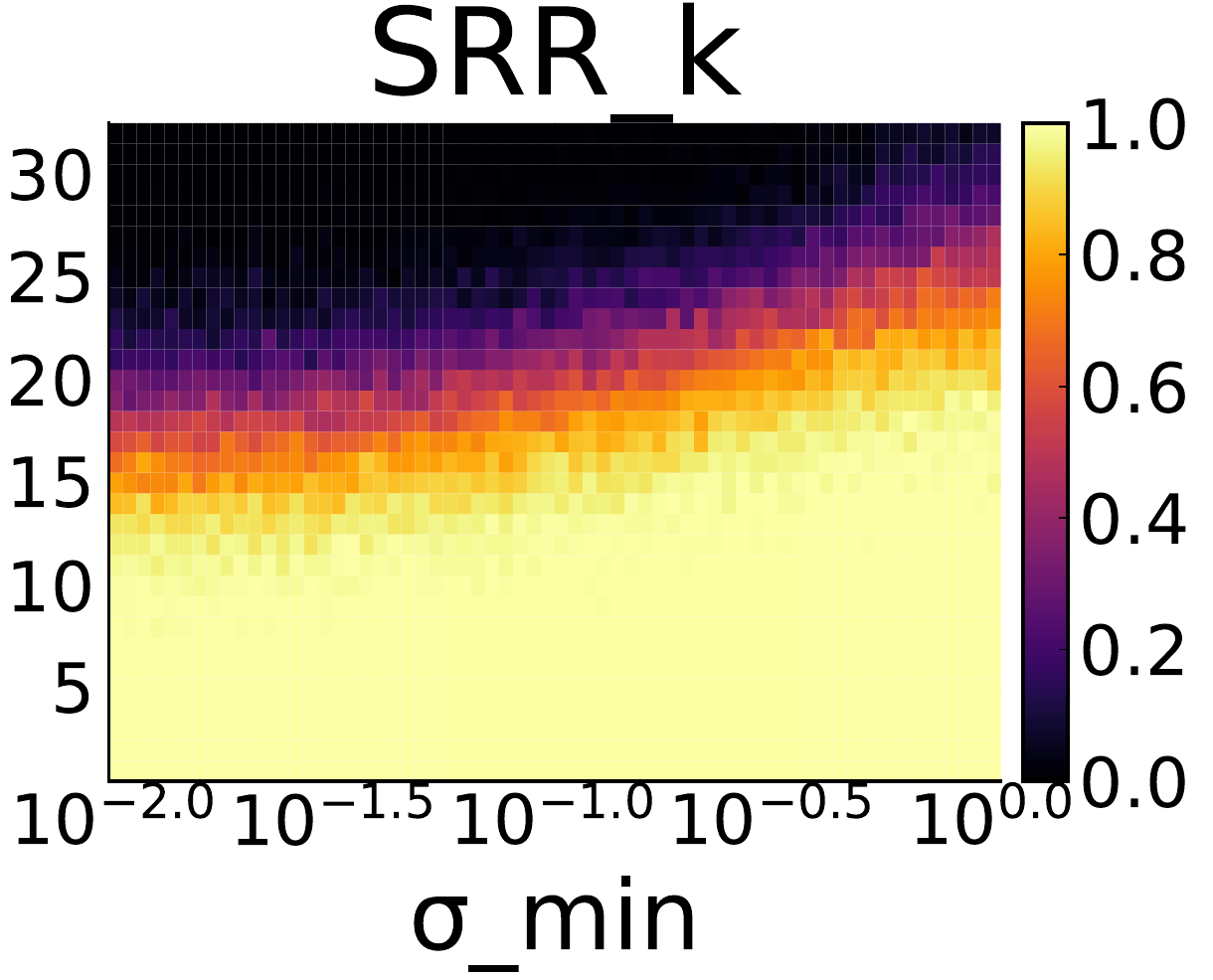}
\caption{Empirical frequency of support recovery as a function of $\sigma_{\min}$ and sparsity $k$ for 64 by 128 matrices.
}
\end{center}
\label{fig:twostage}
\end{figure}

Figure \ref{fig:twostage} shows the empirical frequency of support recovery for the two-stage
algorithms
as a function of sparsity level and the smallest singular value $\sigma_{\min}$ for 64 by 128 matrices.
We generated the problems in a similar fashion as for Figure~\ref{fig:backward}.
OMPR and SRR refer to $s=1$, while SRR$_k$ has $s= k$. 
Among the single step methods,
SRR maintains much better performance on 
matrices with smaller $\sigma_{\min}$ than OMPR.
Further, SRR is marginally less performant than SP,
since SP adds and deletes $k$ atoms,
but is much more efficient than SP,
as it requires far fewer updates - $2$ versus $2k$ - to the QR factorization, which is the most expensive operation.
The most performant method is SRR$_k$ with a small but statistically significant edge over SP for matrices with smaller $\sigma_{\min}$.
Since SP and SRR require a similar effort to implement,
the SRR family provides a performant and elegant alternative to SP
with a potential for higher computational efficiency by controlling the size of the support updates.

\section{CONCLUSION}
We derived novel guarantees for Backward Regression,
which are efficiently computable, 
insightful, 
and experimentally verifiable.
We put forth a new implementation of the backward algorithm
that is both efficient and numerically stable.
Our analysis of LACE
connects the results to existing algorithms, like Subspace Pursuit,
and is also related to magnitude pruning, a principle used to sparsify neural networks. 
Inspired by our theoretical results,
we proposed a new family of two-stage algorithms for compressed sensing problems:
Stepwise Regression with Replacement (SRR).
Remarkably, SRR maintains high performance on coherent dictionaries, is efficient,
and simple to implement.
Lastly, the results herein can help guide users of Backward Regression in professional statistical software, like SAS,
and inspire new research and applications.

\small

\bibliographystyle{IEEEbib}
\bibliography{BackwardOptimalityICASSP}

\end{document}